\documentclass{amsart}

\usepackage{amsmath,amsfonts}

\usepackage{amscd,amsthm}

\usepackage[dvips]{epsfig}

\newcommand{\R}{\mathop{\mathbb{R}}}

\newcommand{\N}{\mathop{\mathbb{N}}}

\theoremstyle{plain} %% This is the default
\newtheorem{theorem}{Theorem}[section]

\newtheorem{remark}[theorem]{Remark}

\numberwithin{equation}{section}

\begin{document}

\title[Qualitative analysis of phase--portrait.]{Qualitative analysis of 
phase--portrait for a class of planar vector fields via the comparison method.}

\author{Timoteo Carletti, Lilia Rosati and Gabriele Villari}

\address[Timoteo Carletti]{Scuola Normale Superiore, piazza dei Cavalieri 7,
 56126  Pisa, Italy}

\address[Lilia Rosati]{Dip. Matematica ''U. Dini'', viale Morgagni 67/A,
 50134 Firenze, Italy}

\address[Gabriele Villari]{Dip. Matematica ''U. Dini'', viale Morgagni 67/A,
 50134 Firenze, Italy}

\email[Timoteo Carletti]{t.carletti@sns.it}

\email[Lilia Rosati]{lrosati@math.unifi.it}

\email[Gabriele Villari]{villari@math.unifi.it}


\begin{abstract}
The phase--portrait of the second order differential equation: $$\ddot
x+\sum_{l=0}^nf_l(x) \dot x^l=0\, ,$$ is studied. Some results concerning
existence, non--existence and uniqueness of limit cycles are
presented. Among these, a generalization of the classical Massera
uniqueness result is proved.
\end{abstract}
\maketitle

\section{Introduction}
\label{sec:intro}

In this paper we investigate the qualitative behavior of the
differential equation:
\begin{equation}
  \label{eq:sode}
  \ddot x +P(x,\dot x)=0 \, ,
\end{equation}
where $P(x,y)=\sum_{l=0}^n f_l(x)y^l$, for some fixed $n\in \N$, which in 
the phase--plane 
can be rewritten as the differential autonomous system:
\begin{equation}
\label{eq:systsode}
\begin{cases}
\dot x &= y \\
\dot y &= -P(x,y) \, .
\end{cases}
\end{equation}
Such systems have been proposed
by H.I. Freedman and Y. Kuang~\cite{freedmankuang} in the study of a
Gause--type predator--prey model.

If we assume some standard regularity assumptions on the functions
$f_l(x)$, $l~\in~\{0,\dots ,n\}$, for the uniqueness of Cauchy's problem,
for instance let $f_0(x)=g(x)$:
\begin{enumerate}
\item[A1)] $g \in {\it Lip}(\R)$;
\item[A2)] $f_l \in \mathcal{C}(\R)$, $l\in \{ 1,\dots ,n \}$;
\end{enumerate}
and the hypothesis:
\begin{enumerate}
\item[B)] $xg(x)>0$ for all $x\neq 0$;
\end{enumerate}
then the origin is the only singular
point, and trajectories turn clockwise around it.

The problem of existence and uniqueness of limit
cycles for system~\eqref{eq:sode} will be considered. All examples
and applications presented will involve polynomials, because this is,
 from our point of view, the more interesting situation; we nevertheless
 observe that all our results hold in a more general setting than polynomial
functions.

Clearly equation~\eqref{eq:sode} is a particular case of the classical
generalized Li\'enard equation:
\begin{equation*}
\ddot x+f(x,\dot x) \dot x+g(x) =0 \, ,
\end{equation*}
as well as the equivalent system:
\begin{equation*}
\begin{cases}
\dot x &= y \\
\dot y &= -f(x,y)y-g(x) \, .
\end{cases}
\end{equation*}
is a generalization of system~\eqref{eq:systsode}.

We observe that a large number of results, in the huge literature concerning
 these generalized Li\'enard equations, are obtained with some assumptions 
on the growth of $f(x,y)$, for instance:
\begin{equation*}
  f(x,y)>-M \, ,
\end{equation*}
being $M$ some suitable positive constant. This because otherwise
trajectories might be not continuable and the uniform boundedness of the
trajectories is a crucial step in order to apply the classical
Poincar\'e--Bendixon Theory. 

This is not the case of system~\eqref{eq:systsode} where in general we
will consider functions $f_l(x)$, $l\in \{0,\dots,n\}$ to be polynomials, 
hence most of the
classical result cannot be used here.

As far as we know system~\eqref{eq:systsode} has been studied only in
two particular cases (the classical Li\'enard equation):
\begin{equation*}
  \label{eq:caso1}
  P_1(x,y) = g(x)+f_1(x)y \, ,
\end{equation*}
and 
\begin{equation*}
  \label{eq:caso2}
  P_2(x,y) = g(x)+f_1(x)y+f_2(x)y^2 \, ,
\end{equation*}
largely studied by several authors in the last decade (see for
instance~\cite{guidorizzi,jifa2,zhou,jifa4,jifa5}). 

We believe that the main reason for this, is the existence of a well known
transformation, discussed in section~\ref{sec:nonexist},
 bringing equation:
\begin{equation}
  \label{eq:equazione5}
  \ddot x+f_1(x)\dot x+f_2(x)\dot x^2 +g(x) =0 \, ,
\end{equation}
into some Li\'enard system. In this way all classical results of existence
 of limit cycles for the
Li\'enard equation may be used. A second reason, for considering 
systems~\eqref{eq:equazione5}, lies in the fact that
its trajectories are always continuable. Of course this can be obtained via the
transformation into a Li\'enard equation, but even in a direct way,
considering the slope of the trajectories in the phase--plane. 

When one is interested in general cases, these properties are no longer valid.

In this paper we will exploit some geometrical properties of the
phase--portrait of system~\eqref{eq:systsode} and when $P(x,y)=P_2(x,y)$,
 we will not make use of the
transformation of the system into a Li\'enard one, but we will use a slight
different system.

The plan of the paper is the following. In
section~\ref{sec:nonexist} we will study system~\eqref{eq:systsode}, splitting
$P$ into two parts: one exhibiting some suitable symmetries and considering
the other as a ''perturbation'', then some properties will be enlighten using a
comparison method. The main result of this
section is the following Theorem, which gives a necessary condition
for the existence of limit cycles.

\begin{theorem}[Non--existence]
  \label{thm:nonexist}
Let us consider the second order differential equation:
\begin{equation}
  \label{eq:general}
  \ddot x +\sum_{l=1}^N f_{2l-1}(x) \dot x^{2l-1} +
 \sum_{l=1}^M f_{2l}(x) \dot x^{2l} +g(x) = 0 \, , 
\end{equation}
where $N$ and $M$ are positive integer. Assume regularity hypotheses A and B 
to be satisfied by the $f_l(x)$'s and $g(x)$, then
 systems~\eqref{eq:general} has not periodic orbits provided all the
 $f_{2l-1}(x)$'s, $l\in\{1 ,\dots,N\}$, never change sign. 
\end{theorem}

Section~\ref{sec:existence} deal with system~\eqref{eq:equazione5}. As already 
mentioned the study is performed using a new differential system rather
 than transforming
system~\eqref{eq:equazione5} into some Li\'enard equation. The main
result of this section is the following:

\begin{theorem}[Existence]
  \label{thm:exist}
Assume hypotheses A and B to hold, let $F_1(x)=\int_0^x f_1(s) \, ds$ 
, $G(x)=\int_0^x g(s)\, ds$ and assume:
\begin{enumerate}
\item[C)] there exists $\delta>0$ s.t. $f_1(x)<0$ for $|x|<\delta$, but $f_1$ is
not always negative;
\item[D1)] there exists $c>0$ s.t.:
\begin{equation*}
F_1(x)>-c \quad \text{if} \,\, x>0 \quad \text{and} \quad 
F_1(x)<c \quad \text{if} \,\, x<0 \, ;
\end{equation*}
\item[D2)] $\limsup_{x\rightarrow \pm \infty}\left[G(x)\pm F_1(x)\right]=
+ \infty$;
\item[D3)] $\lim_{x\rightarrow \pm\infty}\int_0^x f_2(s) \, ds=
l_{\pm}>-\infty$,
\item[E)] there exists $\Delta^{\prime}>0$ s.t. either:
\begin{enumerate}
\item[E1)] $f_2(x) >0$ for $x\leq -\Delta^{\prime}$,
\item[E2)] $\limsup_{x\rightarrow -\infty}\left[F_1(x)+
\sqrt{\frac{-g(x)}{f_2(x)}}\right]=L_{-}<+\infty$;
\end{enumerate}
or
\begin{enumerate}
\item[E1')] $f_2(x) <0$ for $x\geq \Delta^{\prime}$,
\item[E2')] $\limsup_{x\rightarrow +\infty}\left[F_1(x)-
\sqrt{\frac{-g(x)}{f_2(x)}}\right]=L_{+}>-\infty$.
\end{enumerate}
\end{enumerate}
Then system~\eqref{eq:equazione5} has at least a periodic orbit.
\end{theorem}

Certainly this result can be somehow, directly obtained after
 transforming the initial equation into a Li\'enard's one, but
however we point out 
that in our formulation based on a geometric study of the phase--plane, 
assumptions keep their geometrical meaning, as has been used
 in~\cite{villari82}.
 Moreover in some applications, as for instance Theorem~\ref{thm:existpoly}, they can
be easily verified.

Some examples of limit cycles for the general case, when there is not
global boundedness of solutions, are also presented.

In Section~\ref{sec:massera} we consider the problem of uniqueness of
limit cycles. A classical result due to Massera is generalized as to include the
following equation:
\begin{equation}
  \label{eq:massgen}
  \ddot x +\sum_{l=0}^{N}f_{2l+1}(x)\dot x^{2l+1}+x=0 \, ,
\end{equation}
for some positive integer $N$. The Classical Massera result is recovered when
$N=0$.

\begin{theorem}[Existence and Uniqueness]
  \label{thm:masseragen}
Let $f_{2l+1}(x)$, $l\in \{0,\dots, N\}$, $N \geq 1$,verify assumption A and moreover:
\begin{enumerate}
\item[L1)] there exists $\delta>0$ s.t. $f_1(x)<0$ for $|x|<\delta$;
\item[L2)] $f_{2l+1}(x)\geq 0$, $l\in \{1,\dots, N\}$, for all $x$;
\item[L3)] $f_{2l+1}(x)$,  $l\in \{0,\dots, N\}$, is monotone increasing 
(respectively decreasing) for $x>0$ (respectively $x<0$).
\end{enumerate}
Then~\eqref{eq:massgen} has a unique globally attracting limit cycle.
\end{theorem}

As far as we know there are no other generalizations of classical Massera's
result and the problem of adapting his geometrical ideas to a more general
situation, still remain an open question. For this reason we believe
that this particular result is in some way significant.

\begin{remark}
\label{rem:masseraexist}
Assumptions L1 and L3 for $f_1(x)$ don't ensure the existence of the limit cycle
in the Massera's case ($N=0$); on the contrary our assumptions guarantee 
that~\eqref{eq:massgen} always has a limit cycle. For instance if 
$f_1(x)=(x^2-1)e^{-x^2}$ in Massera's case the limit cycle doesn't 
exist~\footnote{\label{ftnt:nonexis}This claim follows easily studying the 
\label{page:nonexis}Massera case in 
the associate Li\'enard plane, where the flow through circles $\{x^2+y^2=r^2\}$
is $-x\int_0^x f_1(s)\, ds$, which doesn't change sign with this choice of 
$f_1(x)$.}
 but adding a ''perturbation'', small as we want, say 
$f_3(x)=\epsilon x^2/(x^2+1)$ with $\epsilon >0$, 
 we ensure existence, and thus uniqueness of the limit cycle.
\end{remark}

\section{A non--existence result}
\label{sec:nonexist}

The aim of this section is to prove the {\em non--existence} result Theorem~\ref{thm:nonexist}.
Throughout this section regularity hypotheses to guarantee existence
and uniqueness of the Cauchy initial problem, will be assumed; also
sign assumption B will hold. 

We will prove this result by considering firstly the simplest case 
of~\eqref{eq:equazione5} dealing with only two terms, and secondly the
general situation involving several terms.

In the phase--plane equation~\eqref{eq:equazione5} can be rewritten as a
first order differential system:
\begin{equation}
  \label{eq:phaseplane}
  \begin{cases}
    \dot u & = v \\
    \dot v & = -f_1(u)v-f_2(u)v^2-g(u) \, ,
  \end{cases}
\end{equation}
thus using the well known transformation~\cite{zhou}:
\begin{equation}
  \label{eq:transform}
  \begin{cases}
   x &= u \\
   y &= \left( v-\int_0^u f_1(s) exp \Big \{ \int_0^s f_2(r) \, dr \Big
   \} \, ds \right) exp \Big \{ -\int_0^u f_2(r) \, dr \Big \} \, ,
  \end{cases}
\end{equation}
with the rescaling of time:
\begin{equation}
  \label{eq:resctime}
d\tau = exp \Big\{ -\int_0^x f_2(s) \, ds \Big\}\, dt \, ,
\end{equation}
 system~\eqref{eq:phaseplane} can be brought
into a Li\'enard one (still denoting by $\dot x = dx/d\tau$): 
\begin{equation}
  \label{eq:lienard}
  \begin{cases}
    \dot x &= y-\tilde{F}(x) \\
    \dot y &= -\tilde{g}(x) \, ,
  \end{cases}
\end{equation}
where:
\begin{equation*}
  \tilde{F}(x)=\int_0^x f_1(s) exp \Big \{ \int_0^s f_2(r) \, dr \Big
   \}\, ds 
   \quad \text{and} \quad \tilde{g}(x)=g(x)exp \Big \{ 2\int_0^x
   f_2(r) \, dr \Big \} \, .
\end{equation*}
A large number of results are based on determine hypotheses on
$\tilde{F}$ and $\tilde{g}$ to ensure existence, uniqueness or
non--existence of limit cycles for system~\eqref{eq:lienard}, and then
to transport them to the original system. We will adopt a different point of 
view, considering system~\eqref{eq:phaseplane} as a ''perturbation'' of the
system:
\begin{equation*}
  \begin{cases}
    \dot u &= v \\
    \dot v &= -f_2(u)v^2-g(u) \, ,
  \end{cases}
\end{equation*}
and exploiting the symmetry of the last one w.r.t. $v\mapsto -v$.

The ''unperturbed'' system, i.e. with $f_1(x) \equiv 0$, is transformed into
 an Hamiltonian one:
\begin{equation}
  \label{eq:hamilt}
  \begin{cases}
    \dot x &= y \\
    \dot y &= -\tilde{g}(x) \, ,
  \end{cases}
\end{equation}
whose Hamilton function is: $H(x,y)=y^2/2+\tilde{G}(x)$, where
$\tilde{G}(x)=\int_0^x \tilde{g}(s)\, ds$. We remark that the origin is
a {\em local center}. It can be a global one according to the
behavior of $\tilde{G}(x)$ for large $|x|$. 

\begin{remark}
  Restricting to the class of polynomials the origin is a {\em
  global center} if and only if $f_2(x)$ is polynomial of odd degree
  with positive leading coefficient. In all the remaining cases the
  origin is a local center and at least a separatrix appears.
\end{remark}

\begin{center}
  \begin{figure}[ht]
   \makebox{
   \epsfig{file=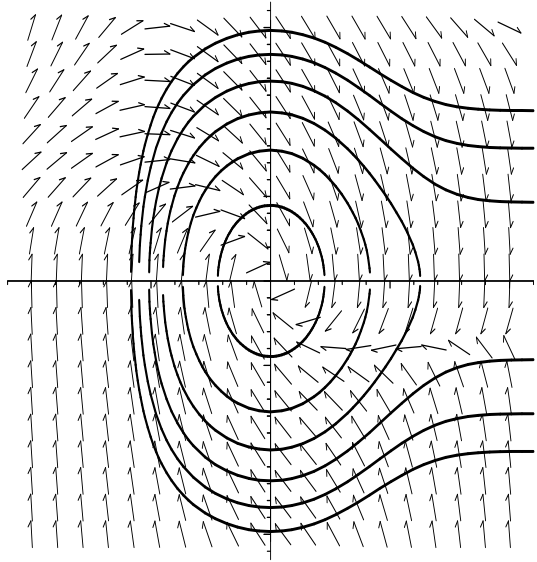,height=7cm,width=7cm}
   }
  \caption{$f_1(x)=x^2+1$, $f_2(x)=-x^2$, $g(x)=x$, energy levels
   $\{ H=\lambda \}$ with $\lambda \in \{ 0.1, 0.3, 0.5, 0.7, 0.9, 1.1\}$.}
  \label{fig:hamdist}
  \end{figure}
\end{center}

The transformation~\eqref{eq:transform} is invertible, hence one can
bring back to the original plane the energy curves $\{ H=const \}$: they will
keep the same topology being only slightly distorted in the
phase--plane (see for example Figure~\ref{fig:hamdist}).

Let us now suppose the perturbation $f_1(x) \not\equiv 0$. We will show, 
using simple
comparison arguments, that system~\eqref{eq:phaseplane} cannot have limit
cycles if $f_1$ never changes sign. Let us suppose $f_1(x) >0$ for all
$x$, and let us compare the {\em slopes} of
system~\eqref{eq:phaseplane}, when $f_1\geq 0$, but not identically
zero, say case A, and $f_1
\equiv 0$ say case B:
\begin{equation*}
  \frac{dv}{du}\Big |_{A}=-f_1(u)+\frac{dv}{du}\Big |_{B}\leq 
  \frac{dv}{du}\Big |_{B} \, . 
\end{equation*}
Hence trajectories in case A enter trajectories of case B and no
limit cycle can exist. This concludes the proof for~\eqref{eq:equazione5}.

 We point out that this situation occurs both
if~\eqref{eq:hamilt} has a global or local center, namely this
situation is not the {\em perturbation} of a global center.

We remark that this is not a new result, in fact if $f_1$ keep constant
sign, then  $\tilde{F}$ vanishes only at $x=0$, and it is positive
(respectively negative) for $x>0$ (respectively $x<0$), thus
Li\'enard system~\eqref{eq:lienard} has not limit cycle and so does
the initial system. But our generalization, obtained from 
Remark~\ref{rem:genertermpari} and Remark~\ref{rem:generdisp}, is to the best
of our knowledge new, and it cannot be obtained ''passing through some
Li\'enard system'' because the previous transformation cannot be
anymore extended.

\begin{remark}
\label{rem:genertermpari}
We can allow terms $f_{2l}(x) \dot x^{2l}$ into~\eqref{eq:equazione5}:
\begin{equation}
\label{eq:general1}
  \ddot x +f_1(x)\dot x+\sum_{l=1}^M f_{2l}(x) \dot x^{2l} + g(x) = 0 \, . 
\end{equation}
For such kind of systems there is not generically a transformation mapping 
them into some Li\'enard systems. Nevertheless we can prove our non--existence
result. A sketch of the proof of this claim 
is as follows: 
assume $f_1\equiv 0$ and consider systems~\eqref{eq:general1} in the
phase--plane. 
 Then there is a trivial symmetry $y\mapsto -y$ which allows to prove the 
existence of a local center. Thus, once we add e perturbation $f_1(x)>0$ we can use the same
reasoning as before, comparing slopes, to prove the non--existence of 
limit cycles.
\end{remark}

The following remark concludes the proof of Theorem~\ref{thm:nonexist}.
\begin{remark}
\label{rem:generdisp}
Provided $f_{2l-1}(x)\geq 0$, $l\in \{ 1,\dots,N\}$ for all $x$, we
can generalize 
further~\eqref{eq:general1} by considering:
\begin{equation*}
  \ddot x +\sum_{l=1}^N f_{2l-1}(x) \dot x^{2l-1} +
 \sum_{l=1}^M f_{2l}(x) \dot x^{2l} +g(x) = 0 \, . 
\end{equation*}
The proof is straightforward and we omit it.
\end{remark}

\section{Existence of limit cycles}
\label{sec:existence}

In this part we will prove two existence results; the first one, 
Theorem~\ref{thm:exist}, is
more general for the involved functions but deals with
system~\eqref{eq:equazione5}, the second one, Theorem~\ref{thm:existpoly},
 holds only for
polynomials but for a slightly more general equation~\eqref {eq:sodepoly}.

\vspace{1em}
\noindent{\bf Theorem~\ref{thm:exist}}\,(Existence){\bf .}
{\em Assume hypotheses A and B to hold, let $F_1(x)=\int_0^x f_1(s) \, ds$ 
, $G(x)=\int_0^x g(s)\, ds$ and assume:
\begin{enumerate}
\item[C)] there exists $\delta>0$ s.t. $f_1(x)<0$ for $|x|<\delta$, but $f_1$ is
not always negative;
\item[D1)] there exists $c>0$ s.t.:
\begin{equation*}
F_1(x)>-c \quad \text{if} \,\, x>0 \quad \text{and} \quad 
F_1(x)<c \quad \text{if} \,\, x<0 \, ;
\end{equation*}
\item[D2)] $\limsup_{x\rightarrow \pm \infty}\left[G(x)\pm F_1(x)\right]=
+ \infty$;
\item[D3)] $\lim_{x\rightarrow \pm\infty}\int_0^x f_2(s) \, ds=
l_{\pm}>-\infty$,
\item[E)] there exists $\Delta^{\prime}>0$ s.t. either:
\begin{enumerate}
\item[E1)] $f_2(x) >0$ for $x\leq -\Delta^{\prime}$,
\item[E2)] $\limsup_{x\rightarrow -\infty}\left[F_1(x)+
\sqrt{\frac{-g(x)}{f_2(x)}}\right]=L_{-}<+\infty$;
\end{enumerate}
or
\begin{enumerate}
\item[E1')] $f_2(x) <0$ for $x\geq \Delta^{\prime}$,
\item[E2')] $\limsup_{x\rightarrow +\infty}\left[F_1(x)-
\sqrt{\frac{-g(x)}{f_2(x)}}\right]=L_{+}>-\infty$.
\end{enumerate}
\end{enumerate}
Then system~\eqref{eq:equazione5} has at least a periodic orbit.}

\begin{remark}[About hypothesis D3]
We point out that hypothesis D3 can be removed, but in that case some growth
conditions on $f_1(x)$ and/or $g(x)$ have to be assumed. The exact ones would be
clear in the proof of Theorem~\ref{thm:exist}.
 Roughly speaking if $\int_0^x f_2(s) \, ds$ diverges to $-\infty$ as some odd
 power, 
$x^{2k+1}$, for $x\rightarrow -\infty$, then $f_1(x)$ has to grow faster than 
$e^{x^{2k+1}}$ to ''compensate'' this divergence. If not the intersection 
property with the $\infty$--isocline is not guarantee, thus there can be
non--winding trajectories.
\end{remark}

\begin{remark}[Polynomial case]
The sign assumptions on $f_2$ (conditions E1 or E1') and hypothesis D3
cannot be simultaneously verified if $f_2(x)$ is a polynomial. Hence 
Theorem~\ref{thm:exist} cannot be applied to the important class of 
polynomials. Nevertheless, we will provide in \S~\ref{ssec:polcase} a result 
(Theorem~\ref{thm:existpoly}) for the following second order
differential equation with polynomial coefficients:
\begin{equation*}
  p(x)\ddot x +p(x)q_1(x)\dot x+q_2(x)\dot x^2 +r(x)=0 \, ,
\end{equation*}
which is strictly related to~\eqref{eq:equazione5}.
\end{remark}

\subsection{Proof of Theorem~\ref{thm:exist}}
\label{ssec:proofthm1}

As already remarked in the previous section, system~\eqref{eq:equazione5}
 rewrites, in the phase--plane:
\begin{equation}
\label{eq:phaseplane2}
  \begin{cases}
    \dot u &= v \\
    \dot v &= -f_{1}(u) v- f_2(u) v^2 -g(u)  \, .
  \end{cases}
\end{equation}

Let us introduce~\cite{papinivillari} a new systems of coordinates
$x=u$ and $y=v+F_1(u)$, such
that the previous phase--plane system rewrites as: 
\begin{equation}
  \label{eq:newsyst}
  \begin{cases}
    \dot x = y -F_1(x) \\
    \dot y = -f_2(x)(y-F_1(x))^2-g(x) \, .
  \end{cases}
\end{equation}

Hypothesis C ensures that the origin is a {\em repellor}, in fact looking at 
the flow through the ovals $\mathcal{O}_r=\{y^2/2+G(x)=r^2\}$ and $r$
is a small positive number, we get: 
\begin{equation*}
\frac{d}{dt}\mathcal{O}_r\Big |_{\it flow}=
-g(x)F_1(x)-yf_2(x)\left(y-F_1(x)\right)^2 \, ,
\end{equation*}
and the claim follows remarking that the last term in the right hand side
is of higher order than $-g(x)F_1(x)$ close enough to the origin, thus the sign
is determined by the first one.

Thus to 
prove the existence of a
closed trajectory is enough to prove the existence of a solution spiraling
 toward the origin.

First of all we want to guarantee that all trajectories intersect the
curve $y=F_1(x)$. For Li\'enard system one has a necessary and
sufficient condition to guarantee such property,
moreover as already remarked in the
previous section, system~\eqref{eq:newsyst} is equivalent, through
the transformation~\eqref{eq:transform} and the rescaling of
time~\eqref{eq:resctime}, to a Li\'enard one : 
\begin{equation}
\label{eq:lienard2}
  \begin{cases}
    \dot \xi &= \eta-\tilde{F}(\xi) \\
    \dot \eta &= -\tilde{g}(\xi) \, ,
  \end{cases}
\end{equation}
where:
\begin{equation*}
  \tilde{F}(\xi)=\int_0^{\xi} f_1(s) exp \Big \{ \int_0^s f_2(r) \, dr \Big
   \}\, ds 
   \quad \text{and} \quad \tilde{g}(\xi)=g(\xi)exp \Big \{ 2\int_0^{\xi}
   f_2(r) \, dr \Big \} \, .
\end{equation*}
We observe that the curve $y=F_1(x)$ corresponds in the 
phase--plane~\eqref{eq:phaseplane2}
to the curve $v=0$, which in the Li\'enard plane~\eqref{eq:lienard2} transform 
into
$\eta=\tilde{F}(\xi)$, hence intersection with $y=F_1(x)$ occurs if
and only if intersection with $\eta=\tilde{F}(\xi)$ occurs.

Let $\tilde{G}(\xi)=\int_0^{\xi} \tilde{g}(s) \, ds$. Assuming
conditions B and D1 one can apply the following
result~\cite{villari}: 
\begin{itemize}
\item for any $\xi_0\geq 0$ and $\eta_0 > \tilde{F}(\xi_0)$, the
  solution of~\eqref{eq:lienard2} passing through $(\xi_0,\eta_0)$
  will intersect $\eta=\tilde{F}(\xi)$ at some
  $\left(\xi^{\prime},\tilde{F}(\xi^{\prime})\right)$, with $\xi^{\prime}>\xi_0$,
  if and only if:
\begin{equation}
  \label{eq:necsuffpos}
\limsup_{\xi \rightarrow +\infty}\left[
\tilde{G}(\xi)+\tilde{F}(\xi)\right] = +\infty \, ;
\end{equation}
\item for any $\xi_0 < 0$ and $\eta_0 < \tilde{F}(\xi_0)$, the
  solution of~\eqref{eq:lienard2} passing through $(\xi_0,\eta_0)$
  will intersect $\eta=\tilde{F}(\xi)$ at some
  $\left(\xi^{\prime},\tilde{F}(\xi^{\prime})\right)$, with $\xi^{\prime}<\xi_0$,
  if and only if:
\begin{equation}
  \label{eq:necsuffneg}
\limsup_{\xi \rightarrow -\infty}\left[
\tilde{G}(\xi)-\tilde{F}(\xi)\right] = +\infty \, .
\end{equation}
\end{itemize}

We claim that from hypotheses D2 and D3 conditions~\eqref{eq:necsuffpos}
and~\eqref{eq:necsuffneg} follow. Let us look for example to $\tilde{G}(\xi)$
 assuming $\int_0^{+\infty} f_2(s)\, ds =l_{+}$ to
be finite and positive, then by definition of limit there exists $M>0$ s.t.:
\begin{equation*}
 \frac{l_{+}}{2}\leq \int_0^x f_2(s)\, ds \leq \frac{3l_{+}}{2} \, ,
\end{equation*}
for all $x\geq M$. Let $\xi > M$, by definition:
\begin{equation*}
\tilde{G}(\xi)=\int_0^M g(x)e^{-2\int_0^x f_2(s)\, ds}\, dx+
\int_M^{\xi} g(x)e^{-2\int_0^x f_2(s)\, ds}\, dx \, ,
\end{equation*}
thus the first term is a constant, depending on $M$, whereas the second can be
estimated by:
\begin{equation*}
e^{-3l_{+}}\int_M^{\xi} g(x)\, dx \leq \int_M^{\xi} g(x)
e^{-2\int_0^x f_2(s)\, ds}\, dx \leq e^{-l_{+}}\int_M^{\xi} g(x)\, dx \, .
\end{equation*}
Finally putting everything together we obtain, for all $\xi > M$, the bound
for $\tilde{G}(\xi)$:
\begin{equation*}
\tilde{G}(M)+e^{-3l_{+}}\left(G(\xi)-G(M)\right) \leq \tilde{G}(\xi) \leq
\tilde{G}(M)+e^{-l_{+}}\left(G(\xi)-G(M)\right)\, ,
\end{equation*}
which implies that $G(\xi)$ and $\tilde{G}(\xi)$ have the same behavior for large
 $\xi$. 

The other cases can be handled similarly and we omit them.

In this way we have the intersection property: trajectory of~\eqref{eq:newsyst}
passing through $(x_0,y_0)$, where $x_0\geq 0$ and $y_0 > F_1(x_0)$, will
intersect the curve $y=F_1(x)$ at some $(x^{\prime},F_1(x^{\prime}))$ with
$x^{\prime} > x_0$. And similarly for the one starting at $(x_0,y_0)$, 
where $x_0< 0$ and $y_0 < F_1(x_0)$.

To conclude the proof we must control that not all trajectories
escape. To do this we will study the
{\em zero--isocline} curves which are explicitly given by:
\begin{equation*}
  y^{\pm}_{\it isoc}(x)=F_1(x)\pm\sqrt{-\frac{g(x)}{f_2(x)}} \, .
\end{equation*}

Let us assume hypothesis E1 to hold (the case for E1' can be handle similarly
and we will omit it). Because of the sign assumptions on $f_2(x)$ and $g(x)$,
$y^{\pm}_{\it isoc}(x)$ are well defined
for $x$ belonging to $(-\infty,-\Delta^{\prime}) \cup
 \mathcal{I}$, where $\mathcal{I}$ is the union of open
disjoint intervals, whose endpoints are zeros~\footnote{If $f_2(0)=0$,
  because also $g(0)=0$, the point $0$ will be an endpoint of some interval in
  $\mathcal{I}$ according to the behavior of $g(x)/f_2(x)$ close to
  $x=0$.} of $f_2(x)$. This gives 
rise to (possible several) branches for the zero--isocline curves. Let
$y^{+}_{left}(x)$ be the leftmost one above $y=F_1(x)$, defined 
for $x<-\Delta^{\prime}$. Such branch can be used, assuming hypothesis E2, to prove 
the existence of a spiraling
 orbit toward the origin.

Before to prove this fact, we have to control the continuability of solutions.
We claim that for all $\bar{\alpha} < \bar{\beta} <0$ (respectively 
$0<\bar{\alpha}^{\prime} < \bar{\beta}^{\prime}$) and for all $\bar{y}>0$ 
(respectively $\bar{y}^{\prime}<0$), the solution of~\eqref{eq:newsyst}
 with initial datum
$(\bar{\alpha},\bar{y})$ will
intersect the line $x=\bar{\beta}$ in finite 
time (respectively the solution of~\eqref{eq:newsyst} with initial datum
$(\bar{\beta}^{\prime},\bar{y}^{\prime})$ will
intersect the line $x=\bar{\alpha}^{\prime}$). To see this just consider 
the system in the 
phase--plane and recall that in this transformation vertical lines are 
mapped on vertical lines. The slope of such a solution is $dv/du =
-f_1(u)-f_2(u)v-g(u)/v$, hence if $|v|$ is large enough, this slope is
sub--linear, giving rise to the continuability from $\bar{\alpha}$ to
$\bar{\beta}$. Moreover, the explicit solution is bounded by some
exponential function which is the responsible for the ''strong''
distortion of the limit cycles shape (see Figure~\ref{fig:exotlc}).

\begin{center}
  \begin{figure}[ht]
   \makebox{
   \epsfig{file=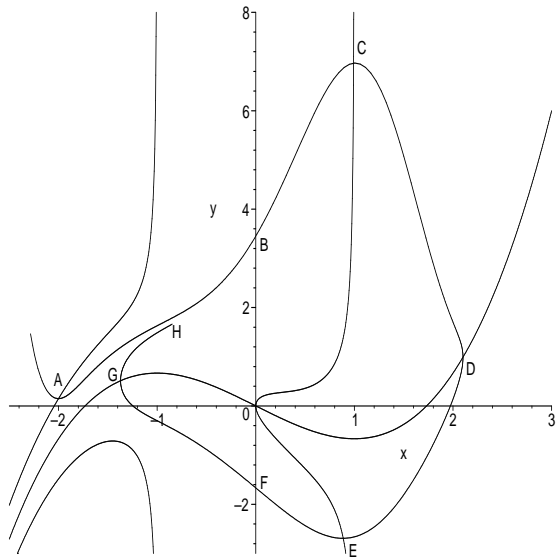,height=9cm,width=9cm,angle=-90}
   }
  \caption{The zero--isocline curves and the $\infty$--isocline curve
   $y=F_1(x)$, 
 with a trajectory spiraling toward the origin.}
  \label{fig:zeroisocl}
  \end{figure}
\end{center}

We are now able to conclude our proof. Let us consider an orbit starting 
at $x_0<-\Delta^{\prime}$, $y_0=y^{+}_{left}(x_0)$, recall that we are assuming 
conditions E1 and E2 to hold, if instead of those we were assuming hypotheses
E1' and E2' we should start from the rightmost branch of the zero--isocline
below $y=F_1(x)$ from some $x_0>\Delta^{\prime}$. Such trajectory will 
intersect the $y$ axis by the previous discussion about continuability (that's
arc $AB$ in Figure~\ref{fig:zeroisocl}). Then, thanks to 
condition~\eqref{eq:necsuffpos} it will reach the curve
$y=F_1(x)$ passing eventually through branches of zero--isocline 
(arcs $BC$ and $CD$ in Figure~\ref{fig:zeroisocl}). On $y=F_1(x)$, $\dot x$ 
changes
 sign and the solution goes back till 
intersecting again the $y$ axis at some negative ordinate
(arcs $DE$ and $EF$ in Figure~\ref{fig:zeroisocl}). Using 
now~\eqref{eq:necsuffneg} we are able to prove that this trajectory will reach 
again the curve $y=F_1(x)$ at some $x<0$, where again $\dot x$ changes sign.
Before reaching this curve, it has to intersect a branch of zero--isocline
below the curve $y=F_1(x)$, hence $\dot y>0$. Because this trajectory
cannot intersect the leftmost branch of $y_{left}^{+}(x)$ 
upon which the vector field points to the right, it will reach anew
the quadrant from which it started spiraling toward
the origin when doing a tour around it. Thus it can be used as outer boundary
of a Poincar\'e--Bendixon domain and this concludes the proof.

\subsection{The polynomial case.}
\label{ssec:polcase}

As already remarked Theorem~\ref{thm:exist} cannot be applied to the class of 
polynomials, in fact the sign assumptions (conditions E1 or E1') 
and hypothesis D3 cannot be simultaneously verified by polynomials. 
Nevertheless we can prove a similar result for a slightly different
second order differential equation with polynomial coefficients.

\begin{theorem}
  \label{thm:existpoly}
The second order differential equation:
\begin{equation}
\label{eq:sodepoly}
  p(x)\ddot x +p(x)q_1(x)\dot x+q_2(x)\dot x^2 +r(x)=0 \, ,
\end{equation}
where $p$, $(q)_{j=1,2}$ and $r$ are {\em real polynomials} of the $x$ variable,
has at least a limit cycle provided the following conditions hold:
\begin{enumerate}
\item[H1)] $p(x)>0$ for all $x$, $q_1(0)<0$ and $xr(x)>0$ for $x\neq 0$;
\item[H2)] $q_1(x)$ has even degree and the leading coefficient is positive;
\item[H3)] $\deg p(x)\geq \deg q_2(x)+2$;
\item[H4)] if $q_2(x)$ has odd degree then the leading coefficient must be
negative;
\item[H5)] $\deg r(x)\leq 2\deg q_1(x)+\deg q_2(x)+1$.
\end{enumerate}
\end{theorem}

We observe that the special form for the coefficient of $\dot x$
is not a restrictive one, we choose it only to state hypotheses of 
Theorem~\ref{thm:existpoly} in a simplest and more clear way than using some
polynomial $\tilde{q}_1(x)$ instead of the factorization $p(x)q_1(x)$.
In this way, the applicability of the present result is straightforward: just
look at the polynomials and compare their degrees.

\proof
Because $p(x)$ never vanishes we can divide~\eqref{eq:sodepoly} by it, then we
 set $f_1(x)=q_1(x)$, $f_2(x)=\frac{q_2(x)}{p(x)}$ and 
$g(x)=\frac{r(x)}{p(x)}$. The proof will be complete once we show that
 $f_1$, $f_2$ and $g$ verify hypotheses of Theorem~\ref{thm:exist}.

Existence and uniqueness of the Cauchy initial value problem associates
 to~\eqref{eq:sodepoly} is trivially guaranteed because we are dealing with 
polynomials.

Signs assumption H1 ensures that hypothesis B holds and by continuity of 
$q_1(x)$ also the first part of assumption C is verified. Also the second part
of it holds, because $q_1$ has even degree and positive leading coefficient.

Hypotheses D1 and D2 are verified because the primitive of $q_1(x)$ is an 
odd degree
polynomial with positive leading coefficient and moreover the primitive of $g(x)$
vanishing at $x=0$ is not negative.

By H3 it follows that the integral $\int_0^x q_2(s)/p(s) \, ds$ is well defined
and finite for all $x$.

Let $\alpha$ (respectively $\beta$) be the smallest (respectively largest)
zero of $q_2(x)$. If $q_2$ has even degree and positive leading coefficient,
 then $q_2(x)$, hence $f_2(x)$, is positive for all $x<\alpha$. 
Moreover let $Q_1(x)=\int_0^x q_1(s) \, ds$, then for $x\rightarrow -\infty$:
\begin{equation*}
Q_1(x)+\sqrt{\frac{-r(x)}{q_2(x)}} \sim x^{\deg q_1 +1}\left(a + 
\sqrt{\frac{b}{c}}|x|^{\frac{\deg r -\deg q_2}{2}-\deg q_1 -1} \right)\, ,
\end{equation*}
where $a$, $b$ and $c$ are the 
positive~\footnote{Assumption on the sign of $r$ implies that $r$ cannot have
even degree nor negative leading coefficient.} leading coefficient of 
respectively $Q_1$, $r$ and $q_2$. Hence assuming H5 hypothesis E2 is verified.

If $q_2$ has even degree and negative leading coefficient,
 then $q_2(x)$, hence $f_2(x)$, is negative for all $x>\beta$. 
Moreover for $x\rightarrow +\infty$ one has:
\begin{equation*}
Q_1(x)-\sqrt{\frac{-r(x)}{q_2(x)}} \sim x^{\deg q_1 +1}\left(a - 
\sqrt{\frac{b}{c}}|x|^{\frac{\deg r -\deg q_2}{2}-\deg q_1 -1} \right)\, ,
\end{equation*}
and again H5 ensures hypothesis E2' to hold.

Finally if $q_2$ has odd degree and negative leading coefficient, then both
couples of hypotheses E1, E2 and E1', E2' hold with 
$\Delta^{\prime}=\max \{|\alpha |,|\beta|\}$.
\endproof

In the following Figure~\ref{fig:exotlc} we show some ''exotic'' limit cycle
whose existence is guaranteed by this last Theorem.

\begin{center}
  \begin{figure}[ht]
   \makebox{
   \epsfig{file=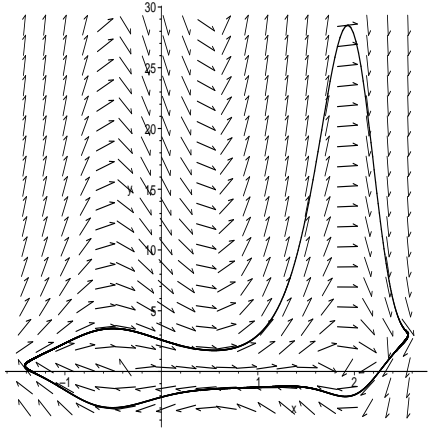,height=7cm,width=7cm,angle=-90}
   }
  \caption{We show an exotic limit cycle, whose existence is proved by 
Theorem~\ref{thm:existpoly}. $p(x)=(x/3)^6+1$, $q_1(x)=x^2-1$, 
$q_2(x)=x^4-4x^2+1$ and $r(x)=x$.}
  \label{fig:exotlc}
  \end{figure}
\end{center}

To conclude this section let us discuss briefly the general
situation where there are few results about existence of limit
cycles. This is because, as already mentioned, one cannot always have
continuability of solutions, hence the 
Poincar\'e--Bendixon Theorem, as
usually done in Literature, cannot be applied ''in
large''. Nevertheless, using standard Hopf's bifurcation argument, we
are able to produce examples with small amplitude limit cycles, even
if global continuability is not allowed.

To fix the ideas, let us consider the following example:
\begin{equation*}
  \ddot x +(ax^2-b)\dot x +(x^2+1)\dot x^2+x^3 \dot x^3+x=0 \, ;
\end{equation*}
phase--space investigation shows that the slope is superlinear, hence
not all trajectories are continuable. Near the origin the dynamics is
dominated by the $\dot x$ term, if $a>0$ and $b=0$ the origin is a
local attractor. For positive, suitable small $b$ a Hopf's bifurcation
takes place: the origin changes its stability and a small amplitude
limit cycle appears, while ''in the large'' orbits behavior doesn't
change too much.

The idea can be generalized as to include other situations, where
always small amplitude limit cycles exist and there is not global
continuability of solutions.

\section{A generalization of a Theorem by Massera}
\label{sec:massera}

\subsection{Introduction}
\label{ssec:masseraintro}

Massera's Theorem~(see \cite{massera} but also~\cite{sansoneconti} pages 307--308)
 concerns the problem of existence and 
uniqueness of limit cycles for the second order differential 
equation~\eqref{eq:massgen}, which in phase--space can
be rewritten as:
\begin{equation}
  \label{eq:massera}
  \begin{cases}
    \dot x &= y \\
    \dot y & =-f_1(x) y -x \, .
  \end{cases}
\end{equation}
Under the assumptions: $f_1(x)<0$ for $|x|<\delta$, for some positive
$\delta$, and $f_1$ increasing (respectively decreasing) for $x>0$
(respectively $x<0$), Massera proved that if a limit cycle exists, then it is
globally attracting, hence unique. The proof is based on the ''geometrical''
construction of the so called {\em geodesic system}, according to Z.--F. Zhang~\cite{zhang}.

\begin{center}
  \begin{figure}[ht]
   \epsfig{file=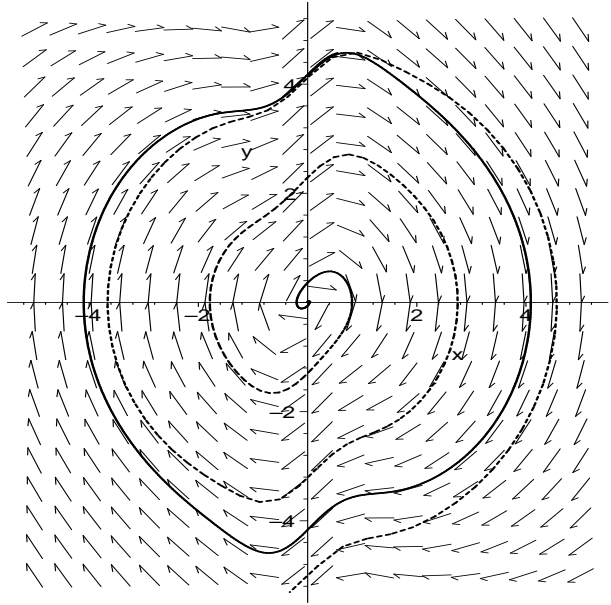,height=8cm,width=8cm,angle=0}
   \caption{The limit cycle for system~\eqref{eq:masseragen} with:
 $f_1(x)=(x^2-1)e^{-x^2}$, 
$f_3(x)=x^2/50(x^2+1)$ (solid line), and the 
spiraling orbit for the classical Massera~\eqref{eq:massera}
 with initial datum $(0,0.01)$ (dashed line).}
   \label{fig:massergen}
  \end{figure}
\end{center}

The previous assumptions on $f_1(x)$ do not guarantee the existence 
 of the limit cycle; an easy calculation 
(see for instance footnote~\ref{ftnt:nonexis} at page~\pageref{page:nonexis}) 
shows that 
system~\eqref{eq:massera} has not limit cycles at all, if for all
$\delta_1<\delta_2$ one has: $\int_{\delta_1}^{\delta_2}f_1(x)\, dx <0$.
We will see in a while that this situation cannot occur in our case.
See Remark~\ref{rem:masseraexist} and Figure~\ref{fig:massergen} for
 an explicit example.

Generalizations of this Theorem, for instance by replacing $\dot y =
-f_1(x)y-x$ with $\dot y =-f_1(x)y-g(x)$, are not available. Here we
propose a new point of view considering the following system:
\begin{equation}
  \label{eq:masseragen}
  \begin{cases}
    \dot x &= y \\
    \dot y & =-f_1(x) y-\sum_{l=1}^N f_{2l+1}(x)y^{2l+1} -x \, ,
  \end{cases}
\end{equation}
under the additional hypothesis: $f_{2l+1}(x)\geq 0$ for all $x$ and
$f_{2l+1}(x)$ is increasing (respectively decreasing) for $x>0$ 
(respectively $x<0$).

We point out that system~\eqref{eq:masseragen} admits a limit cycle even if
system~\eqref{eq:massera} doesn't, moreover if the classical Massera system has 
a unique limit cycle then it must contain in its interior the one 
of~\eqref{eq:masseragen}.

\subsection{Existence of limit cycles}
\label{ssec:existence}

Let us consider the circles $\mathcal{C}_{r}=\{ (x,y)\in \R^2:
x^2+y^2=r^2 \}$, we claim that the flow of~\eqref{eq:masseragen} is
transversal (more precisely it points outward) to $\mathcal{C}_r$ for
$r$ small enough, hence the origin is a {\em repellor}.

To prove this claim let us consider the angle between
$\mathcal{C}_{r}$ and the vector field $X=(y,-f_1(x) y-\sum_{l=1}^N
f_{2l+1}(x)y^{2l+1} -x)$ valued at $\mathcal{C}_{r}$:
\begin{equation*}
\alpha_r= < \nabla \mathcal{C}_r,X> \Big |_{\mathcal{C}_r}=-2y^2\left(
 f_1(x)+\sum_{l=1}^N f_{2l+1}(x)y^{2l} \right)\, ,
\end{equation*}
thus if $r>0$ is small enough, $f_1$ is negative and the remaining
positive terms cannot compensate it.

On the other hand if $r$ is large enough then $\alpha_r<0$, in fact
its leading term, $-2y^{2N+2}f_{2N+1}(x)$, is negative because
of the positiveness of $f_{2N+1}(x)$. Hence by Poincar\'e--Bendixon's 
Theorem we have constructed the inner
and outer boundaries of an invariant domain which must contain at
least one limit cycle.

\subsection{The cycle is star--shaped}
\label{ssec:starshaped}

A fundamental assumption for the proof of Massera's Theorem is that
the limit cycle must be {\em star--shaped} w.r.t. the origin. 
Z.--F. Zhang (see~\cite{zhang} pages 246--253) pointed out, 
that generalizations of Massera's Theorem
following the geometrical ideas of geodetic systems, must rely on this 
assumptions.

The aim of  
this section is to prove that this condition is verified for
system~\eqref{eq:masseragen}. 

\begin{center}
  \begin{figure}[ht]
   \epsfig{file=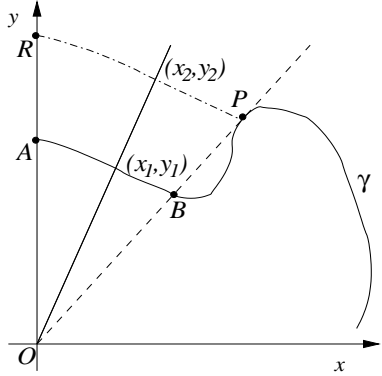,height=5cm,width=5cm}
   \caption{Part of a non star--shaped limit cycle.}
   \label{fig:figura1}
  \end{figure}
\end{center}

Let assume the limit cycle to be not star--shaped w.r.t. the origin,
 thus there is at least a half--line from the 
origin intersecting the cycle at least at two points. Starting from the point $A$
on the $y$--axis the vector field points right--down, $\dot x >0$ and 
$\dot y <0$ if $x$ is small enough, hence the geometry must be as in 
Figure~\ref{fig:figura1}.

Let $\lambda = \overline{PO}/\overline{BO}$ and let us construct the 
arc $RP$ homothetic to $AB$, with ratio $\lambda$. Comparing the slopes of the
vector field on these two arcs we will show a contradiction, hence the limit cycle
must be star--shaped. Take any point $(x_1,y_1)$ on the arc $AB$ and its 
homothetic one $(x_2,y_2)$, let $y=mx$ be the straight line passing through 
these two points. We can assume $x_2 > x_1 >0$. Let us compare the slopes at
$(x_1,mx_1)$ and $(x_2,mx_2)$:
\begin{eqnarray}
\notag 
\frac{dy}{dx}(x_2,mx_2) &=&-f_1(x_2)-\sum_{l=1}^N f_{2l+1}(x_2)(mx_2)^{2l} -m \\
                        &<&-f_1(x_1)-\sum_{l=1}^N f_{2l+1}(x_1)(mx_1)^{2l} -m 
=\frac{dy}{dx}(x_1,mx_1) \, ,
\label{eq:compslopes}
\end{eqnarray}
where we used the fact that the functions 
$x\mapsto -(mx)^{2l}f_{2l+1}(x)$, $l\in \{ 0,\dots,N\}$ are decreasing for
$x>0$. This implies that orbits starting on the arc $RP$ move toward the arc 
$AB$. But this gives rise to a contradiction, in fact for these orbits we have
$\dot x >0$, they cannot accumulate to a point somewhere in the generalized
quadrilateral with vertices $(x_1,y_1)$, $B$, $P$ and $(x_2,y_2)$ because
the origin is the only fixed point, thus they must intersect the limit cycle
 somewhere between points $A$ and $P$, which is impossible because the system
is autonomous and Cauchy's uniqueness and existence assumptions are guaranteed.

\subsection{Uniqueness of limit cycle}
\label{ssec:uniq}

In this paragraph we will prove uniqueness of the limit cycle, by showing its
hyperbolicity. In our proof we will follow closely the original ideas of Massera.

Let us take a limit cycle, $\gamma_1$, and consider its homothetic by a factor 
$\mu >1$. Take a point $(x_1,y_1)\in \gamma_1$ and its homothetic one 
$(x_2,y_2)\in \gamma_2$, let $y=mx$ be the straight line passing through these
two points. Comparing again the slopes, as we did in~\eqref{eq:compslopes}, we get:
\begin{equation*}
\frac{dy}{dx}(x_2,mx_2) \leq \frac{dy}{dx}(x_1,mx_1) \, ,
\end{equation*}
for all $x_1$ and $x_2$. Here the equal sign appears because we are considering
also $x_1=0$. Hence trajectories starting on $\gamma_2$ will go inward 
and accumulating to $\gamma_1$. 

Similarly, considering $\gamma_3$ homothetic of $\gamma_1$ by a factor $0<\mu<1$, and
comparing again slopes we will obtain that $\gamma_1$ is internally attracting.
This proves that $\gamma_1$ is globally hyperbolic, hence unique.

\end{document}